\documentclass[12pt]{article}
\usepackage{amssymb,amsmath}

\newcommand{\Z}{\mathbb{Z}}
\newcommand{\Q}{\mathbb{Q}}
\newcommand{\R}{\mathbb{R}}
\newcommand{\N}{\mathbb{N}}
\newcommand{\F}{\mathcal{F}}
\newcommand{\OO}{\mathcal{O}}
\newcommand{\M}{\mathcal{M}}
\newcommand{\Gal}{\mathrm{Gal}}
\newcommand{\Aut}{\mathrm{Aut}}
\newcommand{\sgn}{\mathrm{sgn}}

\newcommand{\Kbar}{\overline{K}}
\newcommand{\vbar}{\overline{v}}

\newcommand{\mubold}{\boldsymbol{\mu}}
\newcommand{\lambdabold}{\boldsymbol{\lambda}}
\newcommand{\kappabold}{\boldsymbol{\kappa}}
\newcommand{\proof}{\noindent{\em Proof: }}
\newcommand{\qed}{\hspace{\fill}$\square$}
\newcommand{\ra}{\rightarrow}

\DeclareMathOperator{\Char}{char}

\newtheorem{theorem}{Theorem}

\newtheorem{prop}[theorem]{Proposition}

\newtheorem{remark}[theorem]{Remark}
\newtheorem{example}[theorem]{Example}

\numberwithin{equation}{section}
\numberwithin{theorem}{section}

\title{Perturbing Eisenstein polynomials over local
fields}

\author{Kevin Keating \\
Department of Mathematics \\
University of Florida \\
Gainesville, FL 32611 \\
USA \\[.2cm]
{\tt keating@ufl.edu}}

\begin{document}

\maketitle

\begin{abstract}
\noindent
Let $K$ be a local field whose residue field has
characteristic $p$ and let $L/K$ be a finite separable
totally ramified extension.  Let $\pi_L$ be a
uniformizer for $L$ and let $f(X)$ be the minimum
polynomial for $\pi_L$ over $K$.  Suppose
$\tilde{\pi}_L$ is another uniformizer for $L$ such
that $\tilde{\pi}_L\equiv\pi_L+r\pi_L^{\ell+1}
\pmod{\pi_L^{\ell+2}}$ for some $\ell\ge1$ and
$r\in\OO_K$.  Let $\tilde{f}(X)$ be the minimum
polynomial for $\tilde{\pi}_L$ over $K$.  In this paper
we give congruences for the coefficients of
$\tilde{f}(X)$ in terms of $r$ and the coefficients of
$f(X)$.  These congruences improve and extend work of
Krasner \cite{kras}.
\end{abstract}

\section{Introduction} \label{prob}

Let $K$ be a field which is complete with respect to a
discrete valuation $v_K$.  Let $\OO_K$ be the ring of
integers of $K$ and let $\M_K$ be the maximal ideal of
$\OO_K$.  Assume that the residue field
$\Kbar=\OO_K/\M_K$ of $K$ is a perfect field of
characteristic $p$.  Let $K^{sep}$ be a separable
closure of $K$ and let $L/K$ be a finite totally
ramified subextension of $K^{sep}/K$.  Let $\pi_L$ be a
uniformizer for $L$ and let
\[f(X)=X^n-c_1X^{n-1}+\cdots+(-1)^{n-1}c_{n-1}X+(-1)^nc_n\]
be the minimum polynomial of $\pi_L$ over $K$.  Let
$\ell\ge1$, let $r\in\OO_K$, and let $\tilde{\pi}_L$ be
another uniformizer for $L$ such that
$\tilde{\pi}_L\equiv\pi_L+r\pi_L^{\ell+1}
\pmod{\M_L^{\ell+2}}$.  Let
\[\tilde{f}(X)=X^n-\tilde{c}_1X^{n-1}+\cdots
+(-1)^{n-1}\tilde{c}_{n-1}X+(-1)^n\tilde{c}_n\]
be the minimum polynomial of $\tilde{\pi}_L$ over $K$.
In this paper we use the techniques developed in
\cite{sym} to obtain congruences for the coefficients
$\tilde{c}_i$ of $\tilde{f}(X)$ in terms of $r$ and the
coefficients of $f(X)$.

     Let $\phi_{L/K}:\R_{\ge0}\ra\R_{\ge0}$ be the
Hasse-Herbrand function of $L/K$, as defined for
instance in Chapter IV of \cite{cl}.  For $1\le h\le n$
set $k_h=\lceil\phi_{L/K}(\ell)+\frac{h}{n}\rceil$.
Krasner \cite[p.\,157]{kras} showed that for
$1\le h\le n$ we have
$\tilde{c}_h\equiv c_h\pmod{\M_K^{k_h}}$.  In
Theorem~\ref{nochange} we prove that
$\tilde{c}_h\equiv c_h\pmod{\M_K^{k_h'}}$ for certain
integers $k_h'$ such that $k_h'\ge k_h$.  Let $h$ be the
unique integer such that $1\le h\le n$ and $n$ divides
$n\phi_{L/K}(\ell)+h$.  Krasner \cite[p.\,157]{kras}
gave a formula for the congruence class modulo
$\M_K^{k_h+1}$ of $\tilde{c}_h-c_h$.  In
Theorem~\ref{special} we give similar formulas for up to
$\nu+1$ values of $h$, where $\nu=v_p(n)$.

     Heiermann \cite{heier} gave formulas which are
analogous to the results presented here.  Let
$S\subset\OO_K$ be the set of Teichmu\"ller
representatives for $\Kbar$.  Let $\pi_K$ be a
uniformizer for $K$ and let $\F(X)$ be the unique power
series with coefficients in $S$ such that
$\pi_K=\pi_L^n\F(\pi_L)$.  Suppose $\tilde{\pi}_L$ is
another uniformizer for $L$ such that $\tilde{\pi}_L
\equiv\pi_L+r\pi_L^{\ell+1}\pmod{\M_L^{\ell+2}}$ for
some $\ell\ge1$ and $r\in S$.  Let $\tilde{\F}$ be the
series with coefficients in $S$ such that
$\pi_K=\tilde{\pi}_L^n\tilde{\F}(\tilde{\pi}_L)$.  Using
Theorem~4.6 of \cite{heier} one can compute certain
coefficients of $\tilde{\F}$ in terms of $r$ and the
coefficients of $\F$.

     In Section~\ref{digraphs} and we recall some facts
about symmetric polynomials from \cite{sym}.  The main
focus is on expressing monomial symmetric polynomials in
terms of elementary symmetric polynomials.  In
Section~\ref{indices} we define the indices of
inseparability of $L/K$ and some generalizations of the
function $\phi_{L/K}$.  In Section~\ref{pert} we prove
our main results.  In Section~\ref{example} we give some
examples which illustrate how the theorems from
Section~\ref{pert} are applied.

\section{Symmetric polynomials and cycle digraphs}
\label{digraphs}

Let $n\ge1$, let $w\ge1$, and let $\mubold$ be a
partition of $w$.  We view $\mubold$ as a multiset
of positive integers
such that the sum $\Sigma(\mubold)$ of the elements
of $\mubold$ is equal to $w$.  The cardinality of
$\mubold$ is denoted by $|\mubold|$.  For $\mubold$ such
that $|\mubold|\le n$ we let
$m_{\mubold}(X_1,\dots,X_n)$ be the monomial symmetric
polynomial in $n$ variables associated to $\mubold$.
For $1\le h\le n$ let $e_h(X_1,\dots,X_n)$ denote the
elementary symmetric polynomial of degree $h$ in $n$
variables.  By the fundamental theorem of symmetric
polynomials there is a unique polynomial
$\psi_{\mubold}\in\Z[X_1,\dots,X_n]$ such that
$m_{\mubold}=\psi_{\mubold}(e_1,\dots,e_n)$.  In this
section we use a theorem of Kulikauskas and Remmel
\cite{kr} to compute certain coefficients of
$\psi_{\mubold}$.

     The formula of Kulikauskas and Remmel can be
expressed in terms of tilings of a certain type of
digraph.  We say that a directed graph $\Gamma$ is a
cycle digraph if it is a disjoint union of finitely many
directed cycles of length $\ge1$.  We denote the vertex
set of $\Gamma$ by $V(\Gamma)$, and we define the sign
of $\Gamma$ to be $\sgn(\Gamma)=(-1)^{w-c}$, where
$w=|V(\Gamma)|$ and $c$ is the number of cycles that
make up $\Gamma$.

     Let $\Gamma$ be a cycle digraph with $w\ge1$
vertices and let $\lambdabold$ be a partition of $w$.
A $\lambdabold$-tiling of $\Gamma$ is a set $S$ of
subgraphs of $\Gamma$ such that
\begin{enumerate}
\item Each $\gamma\in S$ is a directed path of length
$\ge0$.
\item The collection $\{V(\gamma):\gamma\in S\}$ forms a
partition of the set $V(\Gamma)$.
\item The multiset $\{|V(\gamma)|:\gamma\in S\}$ is
equal to $\lambdabold$.
\end{enumerate}
Let $\mubold$ be another partition of $w$.  A
$(\lambdabold,\mubold)$-tiling of $\Gamma$ is an ordered
pair $(S,T)$, where $S$ is a $\lambdabold$-tiling of
$\Gamma$ and $T$ is a $\mubold$-tiling of $\Gamma$.  Let
$\Gamma'$ be another cycle digraph with $w$ vertices and
let $(S',T')$ be a $(\lambdabold,\mubold)$-tiling of
$\Gamma'$.  An isomorphism from $(\Gamma,S,T)$ to
$(\Gamma',S',T')$ is an isomorphism of digraphs
$\theta:\Gamma\ra\Gamma'$ which carries $S$ onto $S'$
and $T$ onto $T'$.  Say that the
$(\lambdabold,\mubold)$-tilings $(S,T)$ and $(S',T')$ of
$\Gamma$ are isomorphic if there exists an isomorphism
from $(\Gamma,S,T)$ to $(\Gamma,S',T')$.  Say that
$(S,T)$ is an admissible $(\lambdabold,\mubold)$-tiling
of $\Gamma$ if $(\Gamma,S,T)$ has no nontrivial
automorphisms.  Let $\eta_{\lambdabold\mubold}(\Gamma)$
denote the number of isomorphism classes of admissible
$(\lambdabold,\mubold)$-tilings of $\Gamma$.

     Let $w\ge1$ and let $\lambdabold,\mubold$ be
partitions of $w$.  Set
\begin{equation} \label{dlm}
d_{\lambdabold\mubold}
=(-1)^{|\lambdabold|+|\mubold|}\cdot\sum_{\Gamma}\,
\sgn(\Gamma)\eta_{\lambdabold\mubold}(\Gamma),
\end{equation}
where the sum is over all isomorphism classes of cycle
digraphs $\Gamma$ with $w$ vertices.  Since 
$\eta_{\mubold\lambdabold}=\eta_{\lambdabold\mubold}$ we
have $d_{\mubold\lambdabold}=d_{\lambdabold\mubold}$.
Kulikauskas and Remmel \cite[Th.\,1(ii)]{kr} proved the
following:

\begin{theorem} \label{krtheorem}
Let $n\ge1$, let $w\ge1$, and let $\mubold$ be a
partition of $w$ with at most $n$ parts.  Let
$\psi_{\mubold}$ be the unique element of
$\Z[X_1,\dots,X_n]$ such that
$m_{\mubold}=\psi_{\mubold}(e_1,\dots,e_n)$.  Then
\[\psi_{\mubold}(X_1,\dots,X_n)=\sum_{\lambdabold}
d_{\lambdabold\mubold}\cdot
X_{\lambda_1}X_{\lambda_2}\dots X_{\lambda_k},\]
where the sum is over all partitions
$\lambdabold=\{\lambda_1,\dots,\lambda_k\}$ of $w$ such
that $1\le\lambda_i\le n$ for $1\le i\le k$.
\end{theorem}

     We now recall some formulas from \cite{sym} for
computing values of $\eta_{\lambdabold\mubold}(\Gamma)$.

\begin{prop} \label{cycles}
Let $a,b,c,d,w$ be positive integers such
that $a\not=c$, $b\not=d$, and let $r,s$ be
nonnegative integers.  Let $\Gamma$ be a directed cycle
of length $w$. \\[.2cm]
(a) Suppose $w=ra=sb+d$.  Let $\lambdabold$ be the
partition of $w$ consisting of $r$ copies of $a$, and
let $\mubold$ be the partition of $w$ consisting of
$s$ copies of $b$ and one copy of $d$.  Then
$\eta_{\lambdabold\mubold}(\Gamma)=a$.  \\[.2cm]
(b) Suppose $w=ra+c=sb+d$.  Let $\lambdabold$ be the
partition of $w$ consisting of $r$ copies of $a$ and one
copy of $c$, and let $\mubold$ be the partition of $w$
consisting of $s$ copies of $b$ and one copy of $d$.
Then $\eta_{\lambdabold\mubold}(\Gamma)=w$.
\end{prop}

\proof Statement (a) follows from Proposition~2.5 of
\cite{sym} if $s=0$, and from Proposition~2.3 of
\cite{sym} if $s\ge1$.  Statement (b) follows from
Proposition~2.2 of \cite{sym}. \qed \medskip

     Using these formulas we can compute
$d_{\lambdabold\mubold}$ in some cases.

\begin{prop} \label{akcbld}
Let $a,b,c,d,w$ be positive integers such that $a\not=c$
and $b\not=d$.  Let $r,s$ be nonnegative integers such
that $w=ra+c=sb+d$ and $a>sb$.
Let $\lambdabold$ be the partition of
$w$ consisting of $r$ copies of $a$ and 1 copy of $c$,
and let $\mubold$ be the partition of $w$ consisting of
$s$ copies of $b$ and 1 copy of $d$.  Then
\[d_{\lambdabold\mubold}=
\begin{cases}
(-1)^{r+s+w+1}w
&\text{if }b\nmid c\text{ or }sb<c, \\
(-1)^{r+s+w+1}(w-ab)
&\text{if $b\mid c$ and }sb\ge c.
\end{cases}\]
\end{prop}

\proof Let $\Gamma$ be a cycle digraph which has an
admissible $(\lambdabold,\mubold)$-tiling.  Suppose
$\Gamma$ consists of a single cycle of length $w$.
Then by Proposition~\ref{cycles}(b) we have
$\eta_{\lambdabold\mubold}(\Gamma)=w$.  Suppose
$\Gamma$ has more than one cycle.  Since $\Gamma$ has a
$\mubold$-tiling, $\Gamma$ has a cycle $\Gamma_1$ such
that $|V(\Gamma_1)|\le sb$.  Since $a>sb$ and
$\Gamma$ has a $\lambdabold$-tiling, it follows that
$|V(\Gamma_1)|=c=mb$ for some $m$ such that
$1\le m\le s$.  Hence if $\Gamma$ has more than one
cycle we must have $b\mid c$  and $c\le sb$.  Let
$\lambdabold_1$ be the partition
of $c$ consisting of one copy of $c$ and let $\mubold_1$
be the partition of $c$ consisting of $m$ copies of
$b$.  Then every $\lambdabold$-tiling of $\Gamma$
restricts to a $\lambdabold_1$-tiling of $\Gamma_1$, and
every $\mubold$-tiling  of $\Gamma$ restricts to a
$\mubold_1$-tiling of $\Gamma_1$.  It follows from
Proposition~\ref{cycles}(a) that
$\eta_{\lambdabold_1\mubold_1}(\Gamma_1)=b$.

     Let $\Gamma_2$ be another cycle of $\Gamma$.
Since $\Gamma$ has a $\lambdabold$-tiling,
$|V(\Gamma_2)|\ge a>sb$.  Hence every
$\mubold$-tiling of $\Gamma$ restricts to a tiling of
$\Gamma_2$ which includes a path $\delta$ with
$|V(\delta)|=d$.  Since $\mubold$ has only one part
equal to $d$, it follows that
$\Gamma=\Gamma_1\cup\Gamma_2$.  Therefore we have
$|V(\Gamma_2)|=ra=(s-m)b+d$.  Let $\lambdabold_2$ be
the partition of $ra$ consisting of
$r$ copies of $a$ and let $\mubold_2$ be the partition
of $(s-m)b+d=ra$ consisting of $s-m$ copies of $b$
and 1 copy of $d$.  Then every $\lambdabold$-tiling of
$\Gamma$ restricts to a $\lambdabold_2$-tiling of
$\Gamma_2$, and every $\mubold$-tiling  of $\Gamma$
restricts to a $\mubold_2$-tiling of $\Gamma_2$.  It
follows from Proposition~\ref{cycles}(a) that
$\eta_{\lambdabold_2\mubold_2}(\Gamma_2)=a$.  Hence
\[\eta_{\lambdabold\mubold}(\Gamma)=
\eta_{\lambdabold_1\mubold_1}(\Gamma_1)\cdot
\eta_{\lambdabold_2\mubold_2}(\Gamma_2)=ba.\]

     Suppose $b\nmid c$ or $c>sb$.  Then it follows
from the above that the only cycle digraph which has a
$(\lambdabold,\mubold)$-tiling consists of a single
cycle of length $w$.  Hence by (\ref{dlm}) we get
\[d_{\lambdabold\mubold}
=(-1)^{(r+1)+(s+1)}\cdot(-1)^{w-1}w.\]
Suppose $b\mid c$ and $sb\ge c$.  Then $c=mb$ with
$1\le m\le s$.  Hence there are two cycle digraphs
which have a $(\lambdabold,\mubold)$-tiling: a single
cycle of length $w$, and the union of two cycles with
lengths $c=mb$ and $ra=(s-m)b+d$.  Therefore by
(\ref{dlm}) we get
\[d_{\lambdabold\mubold}=
(-1)^{(r+1)+(s+1)}((-1)^{w-1}w+(-1)^{w-2}ab).\]
Hence the formula for $d_{\lambdabold\mubold}$ given in
the theorem holds in both cases. \qed \medskip

     We recall some results from \cite{sym} regarding
the $p$-adic properties of the coefficients
$d_{\lambdabold\mubold}$.  Let $w\ge1$ and let
$\lambdabold$ be a partition of $w$.  For
$k\ge1$ let $k*\lambdabold$ be the partition of $kw$
which is the multiset sum of $k$ copies of
$\lambdabold$, and let $k\cdot\lambdabold$ be the
partition of $kw$ obtained by multiplying the parts of
$\lambdabold$ by $k$.  

\begin{prop} \label{val}
Let $t\ge j\ge0$, let $w'\ge1$, and set $w=w'p^t$.
Let $\lambdabold'$ be a partition of $w'$ and set
$\lambdabold=p^t\cdot\lambdabold'$.  Let $\mubold$
be a partition of $w$ such that there does not exist a
partition $\mubold'$ with $\mubold=p^{j+1}*\mubold'$.
Then $p^{t-j}$ divides $d_{\lambdabold\mubold}$.
\end{prop}

\proof This is proved in Corollary~3.4 of \cite{sym}.
\qed

\begin{prop} \label{cong}
Let $w'\ge1$, $j\ge1$, and $t\ge0$.  Let
$\lambdabold'$, $\mubold'$ be partitions of $w'$ such
that the parts of $\lambdabold'$ are all divisible by
$p^t$.  Set $w=w'p^j$, so that
$\lambdabold=p^j\cdot\lambdabold'$ and
$\mubold=p^j*\mubold'$ are partitions of $w$.  Then
$d_{\lambdabold\mubold} \equiv d_{\lambdabold'\mubold'}
\pmod{p^{t+1}}$.
\end{prop}

\proof This is proved in Proposition~3.5 of \cite{sym}.
\qed

\section{Indices of inseparability} \label{indices}

Let $L/K$ be a totally ramified extension of degree
$n=up^{\nu}$, with $p\nmid u$.  Let $\pi_L$ be a
uniformizer for $L$ whose minimum polynomial over $K$ is
\[f(X)=X^n-c_1X^{n-1}+\cdots+(-1)^{n-1}c_{n-1}X+(-1)^nc_n.\]
For $k\in\Z$ define $\vbar_p(k)=\min\{v_p(k),\nu\}$.
For $0\le j\le\nu$ set
\begin{align} \label{ijpi}
i_j^{\pi_L}
&=\min\{nv_K(c_h)-h:1\le h\le n,\;\vbar_p(h)\le j\} \\
&=\min\{v_L(c_h\pi_L^{n-h}):
1\le h\le n,\;\vbar_p(h)\le j\}-n. \nonumber
\end{align}
Then $i_j^{\pi_L}$ is either a nonnegative integer or
$\infty$; if $\Char(K)=p$ then $i_j^{\pi_L}$ must be
finite, since $L/K$ is separable.  Let $e_L=v_L(p)$
denote the absolute ramification index of $L$.  We
define the $j$th index of inseparability of $L/K$ to be
\begin{equation} \label{ij}
i_j=\min\{i_{j'}^{\pi_L}+(j'-j)e_L:j\le j'\le\nu\}.
\end{equation}
By Proposition~3.12 and Theorem~7.1 of \cite{heier},
$i_j$ does not depend on the choice of $\pi_L$.
Furthermore, our definition of $i_j$ agrees with
Definition~7.3 in \cite{heier}; for the
characteristic-$p$ case see also
\cite[pp.\,232--233]{fried} and \cite[\S2]{fm}.
Write $i_j=A_jn-b_j$ with $1\le b_j\le n$.

\begin{remark} \normalfont
If $i_j^{\pi_L}$ is finite we can write
$i_j^{\pi_L}=a_jn-b_j$ with $a_j\ge1$ (see Section~4
of \cite{sym}).  Thus if $i_j=i_{j'}^{\pi_L}+(j'-j)e_L$
then $A_j=a_{j'}+(j'-j)e_K$.
\end{remark}

     The following facts are easy consequences of the
definitions:
\begin{enumerate}
\item $0=i_{\nu}<i_{\nu-1}\le\dots\le i_1\le i_0<\infty$.
\item If $\Char(K)=p$ then $i_j=i_j^{\pi_L}$.
\item Let $m=\vbar_p(i_j)$.  If $m\le j$ then
$i_j=i_m=i_j^{\pi_L}=i_m^{\pi_L}$.  If $m>j$ then
$\Char(K)=0$ and $i_j=i_m^{\pi_L}+(m-j)e_L$.
\end{enumerate}

     Following \cite[(4.4)]{heier}, for $0\le j\le\nu$
we define functions
$\tilde{\phi}_j:[0,\infty)\ra[0,\infty)$ by
$\tilde{\phi}_j(x)=i_j+p^jx$.  The generalized
Hasse-Herbrand functions
$\phi_j:[0,\infty)\ra[0,\infty)$ are then defined by
\begin{equation} \label{phik}
\phi_j(x)
=\min\{\tilde{\phi}_{j_0}(x):0\le j_0\le j\}.
\end{equation}
Hence we have $\phi_j(x)\le\phi_{j'}(x)$ for
$0\le j'\le j$.  Let
$\phi_{L/K}:[0,\infty)\ra[0,\infty)$ be the usual
Hasse-Herbrand function.  Then by Corollary~6.11 of
\cite{heier} we have $\phi_{\nu}(x)=n\phi_{L/K}(x)$.

     For a partition $\lambdabold=
\{\lambda_1,\dots,\lambda_k\}$ whose parts satisfy
$1\le\lambda_i\le n$ define
$c_{\lambdabold}=c_{\lambda_1}c_{\lambda_2}\dots
c_{\lambda_k}$.  The following is proved in
Proposition~4.2 of \cite{sym}.

\begin{prop} \label{Sigma}
Let  $w\ge1$ and let
$\lambdabold=\{\lambda_1,\dots,\lambda_k\}$ be a
partition of $w$ whose parts satisfy
$1\le\lambda_i\le n$.  Choose $q$ to minimize
$\vbar_p(\lambda_q)$ and set $t=\vbar_p(\lambda_q)$.
Then $v_L(c_{\lambdabold})\ge i_t^{\pi_L}+w$.  If
$v_L(c_{\lambdabold})=i_t^{\pi_L}+w$ and
$i_t^{\pi_L}<\infty$ then $\lambda_q=b_t$ and
$\lambda_i=b_{\nu}=n$ for all $i\not=q$.
\end{prop}

\section{Perturbing $\pi_L$} \label{pert}

In this section we prove our main theorems.  We begin by
applying the results of Section~\ref{digraphs} to the
totally ramified extension $L/K$.  Write
$[L:K]=n=up^{\nu}$ with $p\nmid u$.  Let $\pi_L$,
$\tilde{\pi}_L$ be uniformizers for $L$, with minimum
polynomials over $K$ given by
\begin{align*}
f(X)&=X^n-c_1X^{n-1}+\cdots+(-1)^{n-1}c_{n-1}X+(-1)^nc_n
\\
\tilde{f}(X)&=X^n-\tilde{c}_1X^{n-1}+\cdots
+(-1)^{n-1}\tilde{c}_{n-1}X+(-1)^n\tilde{c}_n.
\end{align*}
Let $1\le h\le n$ and set $j=\vbar_p(h)$.  Define a
function $\rho_h:\N\ra\N$ by
\[\rho_h(\ell)=\left\lceil\frac{\phi_j(\ell)+h}{n}
\right\rceil.\]
Let $\ell\ge1$.  We say $\tilde{f}\sim_{\ell}f$
if $\tilde{c}_h\equiv c_h\pmod{\M_K^{\rho_h(\ell)}}$ for
$1\le h\le n$.  Thus $\sim_{\ell}$ is an equivalence
relation on the set of minimum polynomials over $K$ for
uniformizers of $L$.

     Let $\sigma_1,\dots,\sigma_n$ be the $K$-embeddings
of $L$ into $K^{sep}$.  For each partition $\mubold$
with at most $n$ parts define $M_{\mubold}:L\ra K$ by
\[M_{\mubold}(\alpha)
=m_{\mubold}(\sigma_1(\alpha),\dots,\sigma_n(\alpha)).\]
For $1\le h\le n$ define $E_h:L\ra K$ by
\[E_h(\alpha)
=e_h(\sigma_1(\alpha),\dots,\sigma_n(\alpha)).\]
Then $c_h=E_h(\pi_L)$ and
$\tilde{c}_h=E_h(\tilde{\pi}_L)$.

\begin{prop} \label{spec}
Let $\phi(X)=r_1X+r_2X^2+\cdots$ be a power series
with coefficients in $\OO_K$ such that
$\tilde{\pi}_L=\phi(\pi_L)$.  Then for $1\le h\le n$ we
have
\[E_h(\tilde{\pi}_L)=\sum_{\mubold}
r_{\mu_1}r_{\mu_2}\dots r_{\mu_h}M_{\mubold}(\pi_L),\]
where the sum ranges over all partitions
$\mubold=\{\mu_1,\dots,\mu_h\}$ with $h$ parts.
\end{prop}

\proof This is a special case of Proposition~4.4 in
\cite{sym}. \qed

\begin{prop} \label{setting}
Let $n\ge1$, let $w\ge1$, and let $\mubold$ be a
partition of $w$ with at most $n$ parts.  Then
\[M_{\mubold}(\pi_L)=\sum_{\lambdabold}
d_{\lambdabold\mubold}c_{\lambdabold},\]
where the sum is over all partitions
$\lambdabold=\{\lambda_1,\dots,\lambda_k\}$ of $w$ such
that $1\le\lambda_i\le n$ for $1\le i\le k$.
\end{prop}

\proof This follows from Theorem~\ref{krtheorem} by
setting $X_i=E_i(\pi_L)=c_i$. \qed \medskip

     Let $\ell\ge1$.  Our first main result gives
congruences between the coefficients of $f(X)$ and the
coefficients of $\tilde{f}(X)$ under the assumption
$\tilde{\pi}_L\equiv\pi_L\pmod{\M_L^{\ell+1}}$.

\begin{theorem} \label{nochange}
Let $\pi_L$, $\tilde{\pi}_L$ be uniformizers for $L$ and
let $f(X)$, $\tilde{f}(X)$ be the minimum polynomials
for $\pi_L$, $\tilde{\pi}_L$ over $K$.  Suppose there
are $\ell\ge1$ and $\sigma\in\Aut_K(L)$ such that
$\sigma(\tilde{\pi}_L)\equiv\pi_L\pmod{\M_L^{\ell+1}}$.
Then $\tilde{f}\sim_{\ell}f$.
\end{theorem}

\proof We first show that the theorem holds in the case
where $\tilde{\pi}_L=\pi_L+r\pi_L^{\ell+1}$, with
$r\in\OO_K$.
Let $1\le h\le n$ and set $j=\vbar_p(h)$.  For
$0\le s\le h$ let $\mubold_s$ be the partition of
$\ell s+h$ consisting of $h-s$ copies of 1 and $s$
copies of $\ell+1$.  Then by Proposition~\ref{spec} we
have
\begin{equation} \label{tildech}
\tilde{c}_h=E_h(\tilde{\pi}_L)
=\sum_{s=0}^hM_{\mubold_s}(\pi_L)r^s
=c_h+\sum_{s=1}^hM_{\mubold_s}(\pi_L)r^s.
\end{equation}
To prove that $\tilde{c}_h\equiv c_h
\pmod{\M_K^{\rho_h(\ell)}}$ it's enough to show that
$v_K(M_{\mubold_s}(\pi_L))\ge\rho_h(\ell)$ for
$1\le s\le h$.  Therefore by Proposition~\ref{setting}
it suffices to show
$v_L(d_{\lambdabold\mubold_s}c_{\lambdabold})\ge
\phi_j(\ell)+h$ for all $1\le s\le h$ and all partitions
$\lambdabold$ of $\ell s+h$ whose parts are at most $n$.

     Let $1\le s\le h$, set $j=\vbar_p(h)$, and set
$m=\min\{j,\vbar_p(s)\}$.  Then $m\le j$ and $s\ge p^m$.
Let $\lambdabold=\{\lambda_1,\dots,\lambda_k\}$ be a
partition of $\ell s+h$ such that $1\le\lambda_i\le n$
for $1\le i\le k$.  Choose $q$ to minimize
$\vbar_p(\lambda_q)$ and set $t=\vbar_p(\lambda_q)$.
By Proposition~\ref{Sigma} we have
$v_L(c_{\lambdabold})\ge i_t^{\pi_L}+\ell s+h$.
Suppose $m<t$.  Then $m<\nu$, so we have
$p^{m+1}\nmid\gcd(h-s,s)$.  Hence by
Proposition~\ref{val} we get
$v_p(d_{\lambdabold\mubold_s})\ge t-m$.  Thus
\begin{align*}
v_L(d_{\lambdabold\mubold_s}c_{\lambdabold})&
=v_L(d_{\lambdabold\mubold_s})+v_L(c_{\lambdabold}) \\
&\ge(t-m)v_L(p)+i_t^{\pi_L}+\ell s+h \\
&\ge i_m+\ell p^m+h.
\end{align*}
Suppose $m\ge t$.  Then
\begin{align*}
v_L(d_{\lambdabold\mubold_s}c_{\lambdabold})
&\ge v_L(c_{\lambdabold}) \\
&\ge i_t^{\pi_L}+\ell s+h \\
&\ge i_t+\ell p^m+h \\
&\ge i_m+\ell p^m+h.
\end{align*}
In both cases we get
$v_L(d_{\lambdabold\mubold_s}c_{\lambdabold})
\ge\tilde{\phi}_m(\ell)+h\ge\phi_j(\ell)+h$, and hence
$\tilde{c}_h\equiv c_h\pmod{\M_K^{\rho_h(\ell)}}$.
Since this holds for $1\le h\le n$ we get
$\tilde{f}\sim_{\ell}f$.

     We now prove the general case.  Since $\tilde{f}$
is the minimum polynomial of $\sigma(\tilde{\pi}_L)$
over $K$ we may assume without loss of generality that
$\tilde{\pi}_L\equiv\pi_L\pmod{\M_L^{\ell+1}}$.  By
repeated application of the special case above we get a
sequence
$\pi_L^{(0)}=\pi_L,\pi_L^{(1)},\pi_L^{(2)},\ldots$ of
uniformizers for $L$ with minimum polynomials
$f^{(0)}=f,f^{(1)},f^{(2)},\ldots$ such that for all
$i\ge0$ we have
$\pi_L^{(i)}\equiv\tilde{\pi}_L\pmod{\M_L^{\ell+i+1}}$
and $f^{(i+1)}\sim_{\ell+i}f^{(i)}$.  It follows that
$f^{(i+1)}\sim_{\ell}f^{(i)}$, and hence that
$f^{(i)}\sim_{\ell}f$ for all $i\ge0$.  Since the
sequence $(f^{(i)})$ converges coefficientwise to
$\tilde{f}$ it follows that $\tilde{f}\sim_{\ell}f$.
\qed

\begin{remark} \normalfont
It follows from Theorem~\ref{nochange} that if
$\sigma(\tilde{\pi}_L)\equiv\pi_L\pmod{\M_L^{\ell+1}}$
for some $\sigma\in\Aut_K(L)$ then
$\tilde{c}_h\equiv c_h\pmod{\M_K^{\rho_h(\ell)}}$ for
$1\le h\le n$.  Define functions $\kappa_h:\N\ra\N$ by
\[\kappa_h(\ell)=\left\lceil\frac{\phi_{\nu}(\ell)+h}{n}
\right\rceil.\]
Krasner \cite[p.\,157]{kras} showed that
$\tilde{c}_h\equiv c_h\pmod{\M_K^{\kappa_h(\ell)}}$.
Since $\kappa_h(\ell)\le\rho_h(\ell)$ Krasner's
congruences are in general weaker than the congruences
that follow from Theorem~\ref{nochange}.  However, if
$\ell$ is greater than or equal to the largest lower
ramification break of $L/K$ then
$\phi_j(\ell)=\phi_{\nu}(\ell)$ for $0\le j\le\nu$.
Therefore Theorem~\ref{nochange} does not improve on
\cite{kras} in these cases.
\end{remark}

     For certain values of $h$ we get a more refined
version of the congruences that follow from
Theorem~\ref{nochange}.

\begin{theorem} \label{special}
Let $L/K$ be a finite totally ramified extension of
degree $n=up^{\nu}$.  For $0\le m\le\nu$ write the $m$th
index of inseparability of $L/K$ in the form
$i_m=A_mn-b_m$ with $1\le b_m\le n$.  Let $\pi_L$,
$\tilde{\pi}_L$ be uniformizers for $L$ such that there
are $\ell\ge1$, $r\in\OO_K$, and $\sigma\in\Aut_K(L)$
with $\sigma(\tilde{\pi}_L)\equiv\pi_L+r\pi_L^{\ell+1}
\pmod{\M_L^{\ell+2}}$.  Let $0\le j\le\nu$ satisfy
$\vbar_p(\phi_j(\ell))=j$, and let $h$ be the unique
integer such that $1\le h\le n$ and $n$ divides
$\phi_j(\ell)+h$.  Set $k=(\phi_j(\ell)+h)/n$ and
$h_0=h/p^j$.  Then
\[\tilde{c}_h\equiv c_h+\sum_{m\in S_j}
g_mc_n^{k-A_m}c_{b_m}r^{p^m}\pmod{\M_K^{k+1}},\]
where
\begin{align*}
S_j&=\{m:0\le m\le j,\;\phi_j(\ell)=\tilde{\phi}_m(\ell)\}
\\[.3cm]
g_m&=\begin{cases}
(-1)^{k+\ell+A_m}(h_0p^{j-m}+\ell-up^{\nu-m})
&\text{if }b_m<h \\
(-1)^{k+\ell+A_m}(h_0p^{j-m}+\ell)
&\text{if }h\le b_m<n \\
(-1)^{k+\ell+A_m}up^{\nu-m}&\text{if }b_m=n.
\end{cases}
\end{align*}
\end{theorem}

\proof We first prove that the theorem holds for
$\hat{\pi}_L=\pi_L+r\pi_L^{\ell+1}$.  Let
\[\hat{f}(X)=X^n-\hat{c}_1X^{n-1}+\cdots
+(-1)^{n-1}\hat{c}_{n-1}X+(-1)^n\hat{c}_n\]
be the minimum polynomial for $\hat{\pi}_L$ over $K$.
Let $1\le s\le h$ and let $\lambdabold$ be a partition
of $\ell s+h$ whose parts are at most $n$.  Choose $q$
to minimize $\vbar_p(\lambda_q)$ and set
$t=\vbar_p(\lambda_q)$.  Recall that $\mubold_s$ is the
partition of $\ell s+h$ consisting of $h-s$ copies of 1
and $s$ copies of $\ell+1$.  Since
$\vbar_p(h)=\vbar_p(\phi_j(\ell))=j$ it follows from
the proof of Theorem~\ref{nochange} that
$v_K(d_{\lambdabold\mubold_s}c_{\lambdabold})\ge k$.
Suppose $v_K(d_{\lambdabold\mubold_s}c_{\lambdabold})=k$.
Then the inequalities in the proof of
Theorem~\ref{nochange} must be equalities.  Hence there
is $0\le m\le j$ such that $s=p^m$,
$v_L(c_{\lambdabold})=i_t^{\pi_L}+\ell p^m+h$,
and $\phi_j(\ell)=\tilde{\phi}_m(\ell)$.  In particular,
we have $m\in S_j$.  

     Let $w_m=\ell p^m+h$ and let $\kappabold_m$ be the
partition of $w_m$ consisting of $k-A_m$ copies of $n$
and 1 copy of $b_m$.  By Proposition~\ref{Sigma} we see
that $\lambdabold$ has at most one element not equal to
$n$.  Since $\lambdabold$ is a partition of $w_m$, and
\[w_m=\phi_j(\ell)-i_m+h=(k-A_m)n+b_m,\]
it follows that $\lambdabold=\kappabold_m$.  Hence
$c_{\lambdabold}=c_{\kappabold_m}=c_n^{k-A_m}c_{b_m}$
and $\vbar_p(b_m)=\vbar_p(\lambda_q)=t$.  Using equation
(\ref{tildech}) and Proposition~\ref{setting} we get
\begin{equation} \label{chtilde}
\hat{c}_h\equiv c_h+\sum_{m\in S_j}
d_{\kappabold_m\mubold_{p^m}}c_n^{k-A_m}c_{b_m}r^{p^m}
\pmod{\M_K^{k+1}}.
\end{equation}

     Let $m\in S_j$.  Since
\[j=\vbar_p(\phi_j(\ell))=\vbar_p(\tilde{\phi}_m(\ell))
=\vbar_p(i_m+\ell p^m)\]
and $m\le j$ we get $m\le\vbar_p(i_m)=\vbar_p(b_m)$.
Hence $b_m'=b_m/p^m$ is an integer.  Let $\kappabold_m'$
be the partition of
\[w_m'=(k-A_m)up^{\nu-m}+b_m'=h_0p^{j-m}+\ell\]
consisting of $k-A_m$ copies of $up^{\nu-m}$ and 1 copy
of $b_m'$.  Let $\mubold_{p^m}'$ be the partition of
$w_m'$ consisting of $h_0p^{j-m}-1$ copies of 1 and 1
copy of $\ell+1$.  Since $h\le n$ we have
$up^{\nu-m}>h_0p^{j-m}-1$.  Hence if
$b_m'\not=up^{\nu-m}$ then we can compute
$d_{\kappabold_m'\mubold_{p^m}'}$ using
Proposition~\ref{akcbld}.

     Suppose $b_m<h$.  Then $h_0p^{j-m}-1\ge b_m'$, so
by Proposition~\ref{akcbld} we get
\[d_{\kappabold_m'\mubold_{p^m}'}=(-1)^{k+\ell+A_m}
(h_0p^{j-m}+\ell-up^{\nu-m}).\]
Suppose $h\le b_m<n$.  Then $h_0p^{j-m}-1<b_m'$, so by
Proposition~\ref{akcbld} we get
\[d_{\kappabold_m'\mubold_{p^m}'}=(-1)^{k+\ell+A_m}
(h_0p^{j-m}+\ell).\]
Suppose $b_m=n$, so that $b_m'=up^{\nu-m}$.  Since
$up^{\nu-m}>h_0p^{j-m}-1$, the only cycle digraph which
admits a $(\kappabold_m',\mubold_{p^m}')$-tiling consists
of a single cycle $\Gamma$ of length $w_m'$.  By
Proposition~\ref{cycles}(a) we get
$\eta_{\kappabold_m'\mubold_{p^m}'}(\Gamma)=up^{\nu-m}$.
It then follows from (\ref{dlm}) that
\[d_{\kappabold_m'\mubold_{p^m}'}=(-1)^{k+\ell+A_m}
up^{\nu-m}.\]
Hence in all three cases we have
$d_{\kappabold_m'\mubold_{p^m}'}=g_m$.

     Since $m\le t\le\nu$ it follows from (\ref{ij}) and
(\ref{ijpi}) that
\begin{align}
i_m&\le i_t^{\pi_L}+(t-m)e_L \nonumber \\
nA_m-b_m&\le nv_K(c_{b_m})-b_m+(t-m)e_L \nonumber \\
A_m&\le v_K(c_{b_m})+(t-m)e_K  \nonumber\\
k+1&\le k-A_m+v_K(c_{b_m})+(t-m+1)e_K. \label{last}
\end{align}
Since $p^t\mid b_m$ we have $p^{t-m}\mid b_m'$.
Therefore by Proposition~\ref{cong} we get
\[d_{\kappabold_m\mubold_{p^m}}\equiv
d_{\kappabold_m'\mubold_{p^m}'}\pmod{p^{t-m+1}}.\]
Using (\ref{last}) we see that
\begin{alignat*}{2}
d_{\kappabold_m\mubold_{p^m}}c_n^{k-A_m}c_{b_m}&\equiv
d_{\kappabold_m'\mubold_{p^m}'}c_n^{k-A_m}c_{b_m}
&&\pmod{\M_K^{k+1}} \\
&\equiv g_m c_n^{k-A_m}c_{b_m}&&\pmod{\M_K^{k+1}}.
\end{alignat*}
Therefore the theorem holds when
$\tilde{\pi}_L=\hat{\pi}_L$.

     We now prove the theorem in the general case.  We
may assume that
\[\tilde{\pi}_L\equiv\pi_L+r\pi_L^{\ell+1}
\pmod{\M_L^{\ell+2}}.\]
It follows that
$\tilde{\pi}_L\equiv\hat{\pi}_L\pmod{\M_L^{\ell+2}}$, so
by Theorem~\ref{nochange} we get $\tilde{c}_h\equiv
\hat{c}_h\pmod{\M_K^{\rho_h(\ell+1)}}$.
Since $(\phi_j(\ell)+h)/n=k$ and
$\phi_j(\ell+1)>\phi_j(\ell)$ this implies
$\tilde{c}_h\equiv\hat{c}_h\pmod{\M_K^{k+1}}$.  Hence
the theorem holds for $\tilde{\pi}_L$. \qed

\begin{remark} \normalfont
Suppose $\vbar_p(\phi_j(\ell))=j'\le j$.  Then
$\phi_j(\ell)=\phi_{j'}(\ell)$.  In particular,
$\phi_{\nu}(\ell)=\phi_{j'}(\ell)$ with
$j'=\vbar_p(\phi_{\nu}(\ell))$.  Hence if $1\le h\le n$
and $n$ divides $\phi_{\nu}(\ell)+h$ then
Theorem~\ref{special} gives a congruence for
$\tilde{c}_h$ modulo $\M_K^{k+1}$, where
$k=(\phi_{\nu}(\ell)+h)/n$.  This is the congruence
obtained by Krasner \cite[p.\,157]{kras}.  If $\ell$ is
greater than or equal to the largest lower ramification
break of $L/K$ then $\phi_j(\ell)=\phi_{\nu}(\ell)$ for
$0\le j\le\nu$.  Therefore Theorem~\ref{special} does
not extend \cite{kras} in these cases.
\end{remark}

\section{Some examples} \label{example}

In this section we give two examples related to the
theorems proved in Section~\ref{pert}.  We first apply
these theorems to a 3-adic extension of degree 9:

\begin{example} \normalfont
Let $K$ be a finite extension of the 3-adic field $\Q_3$
such that $v_K(3)\ge2$.  Let
\[f(X)=X^9-c_1X^8+\cdots+c_8X-c_9\]
be an Eisenstein polynomial over $K$ such that
$v_K(c_2)=v_K(c_6)=2$, $v_K(c_h)\ge2$ for $h\in\{1,3\}$,
and $v_K(c_h)\ge3$ for $h\in\{4,5,7,8\}$.  Let $\pi_L$
be a root of $f(X)$.  Then $L=K(\pi_L)$ is a totally
ramified extension of $K$ of degree 9, so we have $u=1$,
$\nu=2$.  It follows from our assumptions about the
valuations of the coefficients of $f(X)$ that the
indices of inseparability of $L/K$ are $i_0=16$,
$i_1=12$, and $i_2=0$.  Therefore $A_0=2$, $A_1=2$,
$A_2=1$, and $b_0=2$, $b_1=6$, $b_2=9$.  We get the
following values for
$\tilde{\phi}_j(\ell)$ and $\phi_j(\ell)$:
\[\begin{array}{|c|c|c|c|c|c|c|}
\hline
\ell\rule{0cm}{.45cm}&\tilde{\phi}_0(\ell)
&\tilde{\phi}_1(\ell)&\tilde{\phi}_2(\ell)
&\phi_0(\ell)&\phi_1(\ell)&\phi_2(\ell) \\
\hline
1&17&15&9&17&15&9 \\
\hline
2&18&18&18&18&18&18 \\
\hline
3&19&21&27&19&19&19 \\
\hline
\end{array}\]

     Now let $\tilde{\pi}_L$ be another uniformizer for
$L$, with minimum polynomial
\[\tilde{f}(X)=X^9-\tilde{c}_1X^8+\cdots+\tilde{c}_8X
-\tilde{c}_9.\]
Suppose $\tilde{\pi}_L\equiv\pi_L\pmod{\M_L^2}$.  Then
by Theorem~\ref{nochange} we get
$\tilde{f}\sim_1f$.  Using the table above we find
that
\begin{alignat*}{3}
\tilde{c}_h&\equiv c_h&&\pmod{\M_K^2}
&&\text{ for }h\in\{1,3,9\}, \\
\tilde{c}_h&\equiv c_h&&\pmod{\M_K^3}
&&\text{ for }h\in\{2,4,5,6,7,8\}.
\end{alignat*}
This is an improvement on \cite{kras}, which gives
$\tilde{c}_h\equiv c_h\pmod{\M_K^2}$ for $1\le h\le9$.
If $\tilde{\pi}_L\equiv\pi_L\pmod{\M_L^3}$ we get
$\tilde{f}\sim_2f$, and hence $\tilde{c}_h\equiv
c_h\pmod{\M_K^3}$ for $1\le h\le9$.  If
$\tilde{\pi}_L\equiv\pi_L\pmod{\M_L^4}$ we get
$\tilde{f}\sim_3f$, and hence
\begin{alignat*}{3}
\tilde{c}_h&\equiv c_h&&\pmod{\M_K^3}
&&\text{ for }1\le h\le8, \\
\tilde{c}_9&\equiv c_9&&\pmod{\M_K^4}.
\end{alignat*}
Since the largest lower ramification break of $L/K$ is
2, the congruences we get for $\ell\ge2$ are the same as
those in \cite{kras}.

     Suppose $\tilde{\pi}_L\equiv\pi_L+r\pi_L^2
\pmod{\M_L^3}$, with $r\in\OO_K$.  By the table above we
get $\vbar_3(\phi_0(1))=0$, $\vbar_3(\phi_1(1))=1$,
$\vbar_3(\phi_2(1))=2$ and $S_0=\{0\}$, $S_1=\{1\}$,
$S_2=\{2\}$.  The corresponding values of $h$ are 1, 3,
9, and we have $h_0=1$, $k=2$ in all three cases.  By
applying Theorem~\ref{special} with $\ell=1$, $j=0,1,2$
we get the following congruences:
\begin{alignat*}{2}
\tilde{c}_1&\equiv c_1+(-1)^{2+1+2}(1+1)c_2r
&&\pmod{\M_K^3} \\
&\equiv c_1-2c_2r&&\pmod{\M_K^3} \\
\tilde{c}_3&\equiv
c_3+(-1)^{2+1+2}(1+1)c_6r^3&&\pmod{\M_K^3} \\
&\equiv c_3-2c_6r^3&&\pmod{\M_K^3} \\
\tilde{c}_9&\equiv c_9+(-1)^{2+1+1}c_9^2r^9
&&\pmod{\M_K^3} \\
&\equiv c_9+c_9^2r^9&&\pmod{\M_K^3}.
\end{alignat*}
Only the congruence for $\tilde{c}_9$ follows from
\cite{kras}.

     Suppose $\tilde{\pi}_L\equiv\pi_L+r\pi_L^3
\pmod{\M_L^4}$.  Then $\vbar_3(\phi_2(2))=2$ and
$S_2=\{0,1,2\}$, which gives $h=9$, $h_0=1$, and $k=3$.
By applying Theorem~\ref{special} with $\ell=2$, $j=2$
we get the following congruence:
\begin{alignat*}{2}
\tilde{c}_9&\equiv c_9+(-1)^{3+2+2}(9+2-9)c_9c_2r \\
&\hspace{2cm}+(-1)^{3+2+2}(3+2-3)c_9c_6r^3
+(-1)^{3+2+1}c_9^2c_9r^9&&\pmod{\M_K^4} \\
&\equiv c_9-2c_2c_9r-2c_6c_9r^3+c_9^3r^9
&&\pmod{\M_K^4}.
\end{alignat*}
Suppose $\tilde{\pi}_L\equiv\pi_L+r\pi_L^4
\pmod{\M_L^5}$.  Then $\vbar_3(\phi_0(3))=0$ and
$S_0=\{0\}$, so we get $h=8$, $h_0=8$, and $k=3$.  By
applying Theorem~\ref{special} with $\ell=3$, $j=0$ we
get the following congruence:
\begin{alignat*}{2}
\tilde{c}_8 &\equiv c_8+(-1)^{3+3+2}(8+3-9)c_9c_2r
&&\pmod{\M_K^4} \\
&\equiv c_8+2c_2c_9r&&\pmod{\M_K^4}.
\end{alignat*}
Again, since the largest lower ramification break of
$L/K$ is 2, the congruences we get for $\ell\ge2$ are
the same as those in \cite{kras}. \qed
\end{example}

     One might hope to prove the following converse to
Theorem~\ref{nochange}: If $\pi_L$, $\tilde{\pi}_L$
are uniformizers for $L$ whose minimum polynomials
satisfy $\tilde{f}\sim_{\ell} f$, then there is
$\sigma\in\Aut_K(L)$ such that
$\sigma(\tilde{\pi}_L)\equiv\pi_L\pmod{\M_L^{\ell+1}}$.
The example below shows that this is not necessarily the
case:

\begin{example} \normalfont
Let $\pi_L$ be a root of the Eisenstein polynomial
$f(X)=X^4+6X^2+4X+2$ over the 2-adic field $\Q_2$.
Then $L=\Q_2(\pi_L)$ is a totally ramified
extension of $\Q_2$ of degree 4, with indices of
inseparability $i_0=5$, $i_1=2$, and $i_2=0$.  We get
the following values for $\tilde{\phi}_j(\ell)$ and
$\phi_j(\ell)$:
\[\begin{array}{|c|c|c|c|c|c|c|}
\hline
\ell\rule{0cm}{.45cm}&\tilde{\phi}_0(\ell)
&\tilde{\phi}_1(\ell)&\tilde{\phi}_2(\ell)
&\phi_0(\ell)&\phi_1(\ell)&\phi_2(\ell) \\
\hline
1&6&4&4&6&4&4 \\
\hline
2&7&6&8&7&6&6\\
\hline
3&8&8&12&8&8&8 \\
\hline
\end{array}\]

     Set $\tilde{\pi}_L=\pi_L+\pi_L^2$, and let the
minimum polynomial for $\tilde{\pi}_L$ over $\Q_2$ be
\[\tilde{f}(X)=X^4-\tilde{c}_1X^3+\tilde{c}_2X^2
-\tilde{c}_3X+\tilde{c}_4.\]
By Theorem~\ref{nochange} we have $\tilde{f}\sim_1 f$,
and hence
\begin{alignat*}{2}
\tilde{c}_1&\equiv0&&\pmod{4} \\
\tilde{c}_2&\equiv6&&\pmod{4} \\
\tilde{c}_3&\equiv-4&&\pmod{8} \\
\tilde{c}_4&\equiv2&&\pmod{4}.
\end{alignat*}
Theorem~\ref{special} gives a refinement of the last
congruence:
\begin{alignat*}{2}
\tilde{c}_4&\equiv2+(-1)^{2+1+1}(2+1-2)\cdot2^{2-1}\cdot6
+(-1)^{2+1+1}\cdot2^{2-1}\cdot2&&\pmod{8} \\
&\equiv2&&\pmod{8}.
\end{alignat*}
Using this refinement we get $\tilde{f}\sim_2f$.

     Using \cite{data} (see also Table~4.2 in \cite{jr})
we obtain a list of the degree-4 extensions of $\Q_2$.
Using the data in this list we find that $L/\Q_2$ is not
Galois, and the only quadratic subextension of $L/\Q_2$
is $M/\Q_2$, where $M=\Q_2(\sqrt{-1})$.  Hence
$\Aut_{\Q_2}(L)=\Gal(L/M)$.  Since the lower
ramification breaks of $L/\Q_2$ are 1, 3, and the lower
ramification break of $M/\Q_2$ is 1, the lower
ramification break of $L/M$ is 3.  Hence if
$\sigma\in\Aut_{\Q_2}(L)$ then
$\sigma(\tilde{\pi}_L)\equiv\tilde{\pi}_L\pmod{\M_L^4}$.
Since $\tilde{\pi}_L=\pi_L+\pi_L^2$ we get
$\sigma(\tilde{\pi}_L)\not\equiv\pi_L\pmod{\M_L^3}$.
\qed
\end{example}

\end{document}